\documentclass{amsart}
\usepackage{amssymb}
\usepackage{epsfig}

\title{Non-Crossing Tableaux}

\author{Pavlo Pylyavskyy}

\date{April 13, 2006}

\address{Department of Mathematics, MIT, Cambridge, MA 02141}

\email{pasha (at) math (dot) mit (dot) edu}


\theoremstyle{plain}
\newtheorem{theorem}{Theorem}
\newtheorem{question}{Question}
\newtheorem{proposition}[theorem]{Proposition}
\newtheorem{lemma}[theorem]{Lemma}
\newtheorem{corollary}[theorem]{Corollary}

\theoremstyle{definition}
\newtheorem{definition}[theorem]{Definition}

\theoremstyle{remark}
\newtheorem{remark}[theorem]{Remark}

\def\X{\,\,\lower2pt\hbox{
\begin{picture}(0,0)%
\includegraphics{figX.pstex}%
\end{picture}%
\setlength{\unitlength}{1973sp}%
\begingroup\makeatletter\ifx\SetFigFont\undefined%
\gdef\SetFigFont#1#2#3#4#5{%
  \reset@font\fontsize{#1}{#2pt}%
  \fontfamily{#3}\fontseries{#4}\fontshape{#5}%
  \selectfont}%
\fi\endgroup%
\begin{picture}(324,324)(589,-973)
\end{picture}
}}
\def\noXv{\,\,\lower2pt\hbox{
\begin{picture}(0,0)%
\includegraphics{figNoX.pstex}%
\end{picture}%
\setlength{\unitlength}{1973sp}%
\begingroup\makeatletter\ifx\SetFigFont\undefined%
\gdef\SetFigFont#1#2#3#4#5{%
  \reset@font\fontsize{#1}{#2pt}%
  \fontfamily{#3}\fontseries{#4}\fontshape{#5}%
  \selectfont}%
\fi\endgroup%
\begin{picture}(316,316)(293,-969)
\end{picture}
}}
\def\noXh{\,\,\lower2pt\hbox{
\begin{picture}(0,0)%
\includegraphics{figNoX2.pstex}%
\end{picture}%
\setlength{\unitlength}{1973sp}%
\begingroup\makeatletter\ifx\SetFigFont\undefined%
\gdef\SetFigFont#1#2#3#4#5{%
  \reset@font\fontsize{#1}{#2pt}%
  \fontfamily{#3}\fontseries{#4}\fontshape{#5}%
  \selectfont}%
\fi\endgroup%
\begin{picture}(316,316)(893,-369)
\end{picture}
}}
\def\noXDU{\,\,\lower2pt\hbox{\input{figDU.pstex_t}}}
\def\noXDD{\,\,\lower2pt\hbox{\input{figDD.pstex_t}}}
\def\noXUD{\,\,\lower2pt\hbox{\input{figUD.pstex_t}}}
\def\noXUU{\,\,\lower2pt\hbox{\input{figUU.pstex_t}}}
\def\noXRR{\,\,\lower2pt\hbox{\input{figRR.pstex_t}}}
\def\noXRL{\,\,\lower2pt\hbox{\input{figRL.pstex_t}}}
\def\noXLR{\,\,\lower2pt\hbox{\input{figLR.pstex_t}}}
\def\noXLL{\,\,\lower2pt\hbox{\input{figLL.pstex_t}}}
\def\XRR{\,\,\lower2pt\hbox{\input{figXRR.pstex_t}}}
\def\XRL{\,\,\lower2pt\hbox{\input{figXRL.pstex_t}}}
\def\XLR{\,\,\lower2pt\hbox{\input{figXLR.pstex_t}}}
\def\XLL{\,\,\lower2pt\hbox{\input{figXLL.pstex_t}}}

\def\TL{\mathit{TL}}

\begin{document}

\begin{abstract}
In combinatorics there is a well-known duality between non-nesting and non-crossing objects. In algebra there are many objects which are {\it {standard}}, for example Standard Young Tableaux, Standard Monomials, Standard Bitableaux. We adopt a point of view that these standard objects are really non-nesting, and we find their non-crossing counterparts. 
\end{abstract}

\maketitle

\section{Introduction}

In 1935 Specht \cite{Sp} constructed irreducible representations of $S_n$ as spaces spanned by certain polynomials, with $S_n$ acting by permuting the variables. In his construction, very special role is played by {\it {Standard Young Tableaux}} - a combinatorial object labeling a basis of the irreducible representations. Standard Young Tableaux proved itself to be an extremely useful tool in studying the representation theory of $S_n$. 

Desarmenien, Kung and Rota were studying a characteristic-free approach to invariant theory. They developed a straightening formula which allowed to decompose an arbitrary bideterminant of given content into linear combination of standard bideterminants, see \cite{DKR} for details. The standard bideterminants are labeled by pairs of generalized Standard Young Tableaux. 

The Standard Monomial Theory was discovered by Hodge and Young, and later was greatly developed by Lakshmibai with coauthors, see \cite{LLM} for a survey. Among other things, the theory describes a basis for the coordinate ring of the Grassmanian. The basis elements of the coordinate ring are called Standard Monomials, and the objects labeling Standard Monomials are very similar in spirit to Standard Young Tableaux. 

In his doctoral thesis \cite{Sch} published in 1901 Issai Schur determined the irreducible polynomial representations of $GL_n(\mathbb C)$. These are labeled by partitions $\lambda$ with at most $n$ parts. Moreover, a basis for such a representation can be labeled by Semi-Standard Young Tableaux of shape $\lambda$. In \cite{Gr} Green describes a construction of such a basis using the bideterminants from \cite{DKR}, following the works of Deruyts \cite{De} and Clausen \cite{Cl}.

A common theme of the results mentioned above is the labeling of certain basis by some {\it {standard}} objects. In this paper we view an object a standard if it is {\it {non-nesting}}. In combinatorics there is a well-known duality between non-nesting and non-crossing objects; see for example \cite{CDDSY}. This duality suggests that there might be non-crossing counterparts of the standard bases. We construct such counterparts.

The paper goes as follows. In Section \ref{NCT} we review the construction of Specht modules. Then we define Non-Crossing Tableaux which will be the key object in the paper. Then we study properties of Non-Crossing Tableaux. In particular, we show that they provide a natural labeling for a basis of the irreducible representations of $S_n$. We also show that Semi-Non-crossing Tableaux can be naturally defined, by analogy with Semi-Standard Young Tableaux. The Semi-Noncrossing Tableaux play an important role in further sections. Finally, we show a curious connection with the Temperley-Lieb algebra in the case of two-row Tableaux. In Section \ref{sec:bi} we provide a background on Standard Bitableaux of Rota et.al. Then we define Non-Crossing Bitableaux and show that they possess the basis property. In Section \ref{sec:mon} we start with a brief review of Standard Monomial Theory. We proceed by defining Non-Crossing Monomials and showing that they do form a basis for the coordinate ring of the Grassmannian. Then we pose the question whether Non-Crossing Monomials can be realised as non-initial monomials under some monomial order. We answer the question affirmatively in case of $\mathbb G_{2,n}$. In Section \ref{GL} we describe a basis labeled by NCT for irreducible representations of $GL_n$. Finally, in Section \ref{goodbye} some concluding remarks are made.

\section{Standard and Non-Crossing Tableaux} \label{NCT}
\subsection{Background on Representation Theory of $S_n$}

The following description of Specht's construction is taken from \cite{M}. 

Let $\lambda \vdash n$ be a partition. Recall that a Young diagram is the corresponding shape made out of unit squares. A filling of these squares with numbers from $1$ to $n$ is called a {\it {Tableau}}. If additionally numbers are increasing along rows and columns, then we call such Tableau {\it {a Standard Young Tableau}}, or just SYT.

Let $T$ be a Tableau of shape $\lambda$. Associate to it a polynomial $$P_T = \prod_{i <_T j} (x_i-x_j),$$ where $i <_T j$ means that $i$ is above $j$ in a column of $T$. Let the symmetric group $S_n$ act on the set of such $P_T$ according to the rule $\omega P_T = P_{\omega(T)}$, where $S_n$ acts on Tableaux by permuting the numbers in the filling. Let $S^{\lambda}$ be the $S_n$-module spanned by $P_T$-s, as $T$ runs through all possible tableaux of shape $\lambda$; we denote this set by $\Pi_{\lambda}$.

\begin{theorem}
{\rm [Theorem~1.6.27, Theorem~1.6.29, \cite{M}]{}} \
As $\lambda$ runs through partitions of $n$, $S^{\lambda}$-s form a complete set of irreducible representations of $S_n$. As $T$ runs through SYT-s of shape $\lambda$, $P_T$-s form a basis for $S^{\lambda}$.
\end{theorem}

There exists a way of expressing any $P_T$ in terms of polynomials corresponding to SYT. It is achieved through repeated application of {\it {Garnir relations}} defined as follows. For two columns $a_1 < a_2 < \ldots < a_l$ and $b_1 < b_2 < \ldots < b_m$, find the smallest $r$ such that $a_r > b_r$. Let $H$ be the group of permutations of $a_r, \ldots, a_l, b_1, \ldots, b_r$, and let the $K$ be the subgroup of $H$ which is a product of subgroups permuting $a_r, \ldots, a_l$ and $b_1, \ldots, b_r$. Let $U$ be a system of representatives of cosets of $K$ in $H$.

\begin{theorem}
{\rm [Lemma~1.6.30, \cite{M}]{}} \
$P_T = - \sum_{u \in U} P_{u(T)}$.
\end{theorem}
 
\

We are going to use a variation of the above construction. 

In particular, given a shape $\lambda \vdash n$, consider a rectangular shape $\mu$ which contains $\lambda$. Unless specified otherwise, $\mu$ is going to be the smallest such shape. Now fix a particular filling $F$ of $\mu/\lambda$ with numbers $n+1, \ldots, |\mu|$ in which numbers increase along rows and columns. Thus, any SYT $T$ of shape $\lambda$ can be completed with $F$ to become a SYT $T \cup F$ of shape $\mu$. We denote the set of such resulting SYT-s by $\Pi_{\mu, F}$. Observe that $S_{n}$ via its natural embedding into $S_{|\mu|}$ acts on $\Pi_{\mu, F}$, inducing the action of $S_n$ on the module generated by $P_T$, $T \in \Pi_{\mu, F}$. Denote this module by $S^{\lambda, F}$.

\begin{figure}
\begin{center}
\input{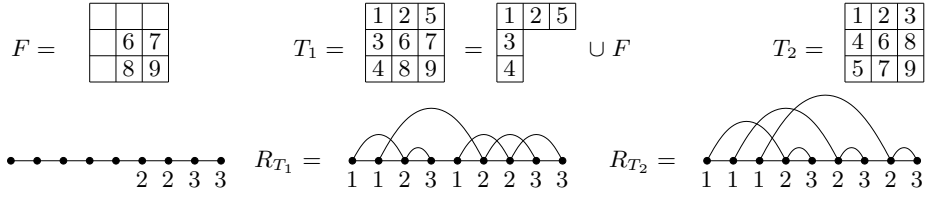}
\end{center}
\caption{$\lambda = (3,1,1)$, $\mu = (3,3,3)$; $T_1$ agrees with $F$; $T_2$ does not.}
\end{figure}

\begin{theorem}
$S_n$-modules $S^{\lambda, F}$ and $S^{\lambda}$ are isomorphic, with an isomorphism given by $P_T \mapsto P_{T \cup F}$.
\end{theorem}

\begin{proof}
One can see that the Garnir relations which hold for $P_T$-s also hold for $P_{T \cup F}$. Therefore, $S^{\lambda, F}$ is a quotient of $S^{\lambda}$. However, since $S^{\lambda}$ is irreducible this implies that they are in fact equal. 
\end{proof}

\begin{remark}
An alternative way to view $S^{\lambda, F}$ is the following: $S^{\lambda, F}$ is the quotient of $S^{\mu}$ by relations $P_T = 0$ for all $T$ which do not agree with $F$.
\end{remark}

\begin{remark}
As was pointed out to the author by Arun Ram, this theorem follows from the results of Garsia and Wachs, \cite{GW}. In their terminology, splitting of $\mu$ into $\lambda$ and $\mu/\lambda$ is an instance of {\it {segmentation}}. It corresponds naturally to the restriction from $S_{|\mu|}$ to the Young subgroup $S_n \times S_{|\mu/\lambda|}$. As a result the irreducible representation of $S_{|\mu|}$ labeled by $\mu$ decomposes into the cross product of the $S_n$-representation $S^{\lambda}$ and the $S_{|\mu/\lambda|}$-representation labeled by the skew shape $\mu/\lambda$. Thus our choice of the particular filling $F$ is just a choice of a particular vector in the letter skew representation.
\end{remark}

As we can see, SYT provide a natural labeling for a basis of $S^{\lambda}$. In the next section we define an alternative combinatorial object which also provides a natural labeling of a (different) basis of $S^{\lambda}$.

\subsection{Non-Crossing Tableaux}

Columns of a Tableau $T$ of rectangular shape $\mu = p \times q$ can be viewed as partitions of numbers from $1$ to $N = pq$ into $q$ sets of size $p$ each. Let us extend the name Tableau to any such partition regardless of particular presentation in terms of filling of some shape. In that case set $\Pi_{\mu, F}$ can be viewed as Tableaux which {\it {agree}} with $F$. By abuse of notation we are sometimes going to say that elements of $\Pi_{\mu, F}$ have shape $\lambda$ (rather than $\mu$); it is going to be clear from the context which one is meant.

Consider two sets of numbers $a$ and $b$ of size $p$ each. Let $a = a_1 < a_2 < \dotsc < a_p$ and $b = b_1 < b_2 < \dotsc < b_p$ be their increasing rearrangements. Call a pair of segments $[x,y]$ and $[z,t]$ {\it {non-crossing}} if following two statements fail to hold: 

\begin{enumerate}
\item $x<z<y<t$;
\item $z<x<t<y$.
\end{enumerate}

Call pair of segments $[x,y]$ and $[z,t]$ {\it {non-nesting}} if the following two statements fail to hold: 

\begin{enumerate}
\item $x<z<t<y$;
\item $z<x<y<t$.
\end{enumerate}

Call $a$ and $b$ non-crossing if for all $1 \leq i \leq p-1$, the segments $[a_i, a_{i+1}]$ and $[b_i, b_{i+1}]$ are non-crossing. Similarly, call $a$ and $b$ non-nesting if for all $1 \leq i \leq p-1$, the segments $[a_i, a_{i+1}]$ and $[b_i, b_{i+1}]$ are non-nesting.

\begin{definition}
A {\it {Non-Crossing Tableau}}, or NCT, is a Tableau in which any two out of $q$ parts are arranged increasingly and are non-crossing.
\end{definition}

\begin{remark}
Formally, NCT are partitions, however the term ``non-crossing partition'' is reserved in the literature for a different notion of ``non-crossing''. 
\end{remark}

Observe, that in this terminology SYT would be called Non-Nesting Tableau. Indeed, the condition which holds for columns of SYT is that any two columns $a$ and $b$ are non-nesting.

\begin{figure}
\begin{center}
\input{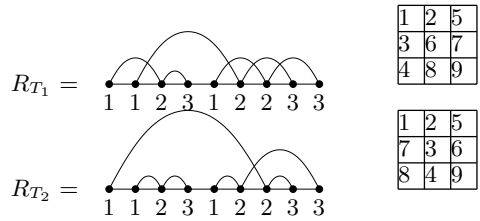}
\end{center}
\caption{$T_1$ is a SYT, $T_2$ is a NCT.}\label{nct5}
\end{figure}

An example of a SYT and NCT can be seen on Figure \ref{nct5}.

It is possible to get rid of $F$ in the definition of NCT.  As before, consider two sets of numbers $a$ and $b$. Let $a = a_1 < a_2 < \dotsc < a_p$ and $b = b_1 < b_2 < \dotsc < b_q$ be their increasing rearrangements, and assume that $p \geq q$. Call pair of segments $[x,y]$ and $[z,t]$ {\it {non-crossing}} if the following two statements fail to hold: 

\begin{enumerate}
\item $x<z<y<t$;
\item $z<x<t<y$.
\end{enumerate}

Call $a$ and $b$ non-crossing if the following conditions hold:

\begin{enumerate}
\item for all $1 \leq i \leq q-1$, segments $[a_i, a_{i+1}]$ and $[b_i, b_{i+1}]$ are non-crossing;
\item if $p>q$ then $a_q < b_q < a_{q+1}$ fails to hold;
\item if $p=q$ then $a_s < b_s$, where $s$ is the largest index such that $a_s \not = b_s$.
\end{enumerate}

Call a Tableau non-crossing if its columns are non-crossing (if we deal with columns of equal length, in the above definition $a$ should correspond to the left one).

\begin{lemma}
For any $F$ the previous definition of NCT is equivalent to the restated definition, up to the order of columns of equal length. 
\end{lemma}

\begin{proof}
One can check that the only influence $F$ has on non-crossing condition is setting an unique order on columns of equal length. Thus if the statement of the lemma is true for some $F$, it is true for any such $F$. Start numbering outer corners of our shape from right to left; when done, do it again for obtained shape, etc. One can check that the previouse definition of a NCT translates into the new definition. 
\end{proof}

It is convenient to look at Tableaux in the following way. Start with a Tableau $T$, whose parts are arranged increasingly. Mark $N$ points on a line (where $N=|\mu|$), which can be identified with integer points on the segment $[1,N]$. To the point $i$ associate a label $l(i)$, such that $i$ is the $l(i)$-th element in the part of $T$ which it belongs to. Call such diagram a {\it {reading}} of $T$, and denote it by $R_T$. For each part $a = a_1 < a_2 < \dotsc < a_p$, connect $a_i$ and $a_{i+1}$ with an arc, $1 \leq i \leq p-1$. Then these arcs define $T$ uniquely.

Here is an example of a NCT with a corresponding reading diagram, for one possible choice of labeling of $F$.

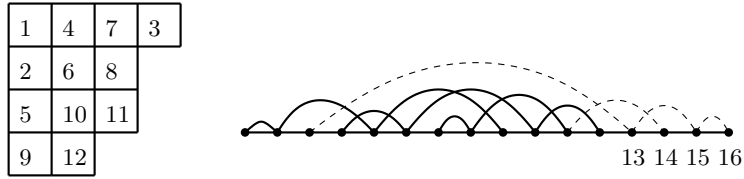
\begin{center}
\input{nct18.pstex_t}
\end{center}

One can see that every reading has the Yamanuchi property: for all $1 \leq k \leq N$ and all $1 \leq m \leq p-1$, the number of points labeled $m$ in $[1,k]$ is not less than the number of points labeled $m+1$ in $[1,k]$. On the other hand, to a particular Yamanuchi word with $q$ $i$-s, $1 \leq i \leq p$ (we consider such words further on, unless specified otherwise), we can always associate a Tableau, possibly in more than one way.

\begin{theorem}
Each (Yamanuchi) reading is a reading of exactly one SYT and exactly one NCT. 
\end{theorem}

\begin{proof}
For each $1 \leq i \leq p$, restrict our attention to points labeled $i$ and $i+1$. There is exactly one way to match them so that the matching is non-crossing. Indeed, put parentheses around neighboring pairs $(i,i+1)$. Then delete those labels that already got matched and repeat. We get a non-crossing matching such that in each pair $i$ precedes $j$, and it is unique.

There is also exactly one way to match them so that the matching is non-nesting. Indeed, match the first $i$ with the first $i+1$, the second $i$ with the second $i+1$, etc. One can see that the resulting matching is non-nesting and no other non-nesting matching exists in which in each matched pair $i$ precedes $i+1$.

In fact, Yamanuchi words in two letters correspond to Dyck paths of length $2q$, and which can be easily bijected with non-crossing and non-nesting matchings.

Now, once we matched the points labelled $i$ and $i+1$ for each $i$ by drawing corresponding arcs, putting those arcs on one picture produces exactly the needed Tableaux. The uniqueness follows from uniqueness for each $i$.
\end{proof}

The illustration of SYT and NCT corresponding to the same reading is given in Figure \ref{nct5}.

\begin{corollary} \label{eq}
The number of NCT of shape $\lambda$ is equal to the number of SYT of shape $\lambda$. 
\end{corollary}

\begin{proof}
Fixing the filling of $F$ is equivalent to fixing the labels of the last $N-n$ points of the reading. Thus, both SYT and NCT of shape $\lambda$ correspond to all possible Yamanuchi words with this particular ending. This correspondance actually provides a bijection from one set of Tableaux to the other.
\end{proof}

\begin{figure}
\begin{center}
\input{nct12.pstex_t}
\end{center}
\caption{}\label{nct12}
\end{figure}

It is a tradition to view Tableaux as fillings of Young diagrams. It is possible to describe a bijection between SYT and NCT in these terms. Namely, for every two consecutive rows of a SYT, write the numbers in those rows in the increasing order on a circle, going clockwise. Mark the numbers from the upper row with a star. Then there is a unique matching of those $2q$ numbers which is non-crossing and such that in each edge the smaller label has a star. The edges of this matching then determine exactly how to arrange the numbers into an NCT. The example is given in Figure \ref{nct12} (Figure \ref{nct5} can be viewed as an example as well). Here of course we could have used lines instead of circles, as on previous figures.

\subsection{Properties of NCT}

By analogy with $\Pi_{\mu, F}$, define $\Xi_{\mu, F}$ to be the set of NCT agreeing with $F$. The elements of the set $\Xi_{\mu, F}$ for $\lambda = (2,1,1)$ and the filling $F$ (unique in this case) are shown in Figure \ref{nct6}. The labels of $F$ are marked. The following theorem shows that NCT-s correspond to a natural basis of $S^{\lambda}$.

\begin{figure}
\begin{center}
\input{nct6.pstex_t}
\end{center}
\caption{}\label{nct6}
\end{figure}

\begin{theorem}
\label{thm:basis}
As $T$ runs through all elements of $\Xi_{\mu, F}$, the corresponding $P_T$-s provide a basis for $S^{\lambda, F}$.
\end{theorem}

\begin{proof}
Since the number of such $P_T$-s is equal to the dimension of the module, it is enough to show that these elements are linearly independent. Introduce the following order on NCT: $T_1 < T_2$ iff $R_{T_1}<R_{T_2}$ in lexicographic order. Suppose $P_T$-s, as $T \in \Xi_{\mu, F}$, are not linearly independent. Then a non-trivial linear combination of these $P_T$-s equals $0$. Choose {\it {the smallest}} (in the order defined above) $T_0$ such that $P_{T_0}$ is involved in one such linear relation with a non-zero coefficient.

Now do the following: assign to each variable $x_i$ the value equal to $l(i)$ - the label of $i$ in $R_{T_0}$. The value of $P_{T_0}$ under assignment is non-zero. Therefore, in the linear dependence relation involving $P_{T_0}$ there should be at least one more $P_T$ which is non-zero under this assignment of values to variables. However, one can see that under this assignment $T_0$ is the largest in the order defined above among the $T$-s such that $P_T$ is nonzero - contradiction.
\end{proof}

Figure \ref{nct7} gives an example of decomposition into $P_T$ corresponding to NCT $T$-s for $\lambda = (2,2,2)$. This corresponds to the equality $$(y-u)(t-u)(y-t)(x-w)(z-w)(x-z) = (x-y)(x-w)(y-w)(z-t)(t-u)(z-u)-$$ $$(x-y)(x-z)(y-z)(t-u)(t-w)(u-w)+(x-t)(x-u)(t-u)(y-z)(z-w)(y-w).$$

\begin{figure}
\begin{center}
\input{nct7.pstex_t}
\end{center}
\caption{}\label{nct7}
\end{figure}

Recall that a Semi-Standard Young Tableau is a Tableau weakly increasing in rows and strictly increasing in columns. One way to look at it is the following: for a given filling $\mu = (\mu_1, \mu_2, \ldots, \mu_k)$, we consider segments $[1, \mu_1], [\mu_1+1, \mu_2], \ldots, [\mu_{k-1}+1, \mu_k]$ of a reading of SYT. Next, we consider all SYT such that each of its columns contains at most one element of each segment. Then we define an equivalence relation on this set where two SYT are equivalent if one is obtained from the other by the action of the Young subgroup $S_{|\mu_1|} \times S_{|\mu_2|} \times \dotsc \times S_{|\mu_k|}$ which permutes the elements of the mentioned segments. Then the equivalence classes are exactly the Semi-Standard Young Tableaux of weight $\mu$.

Remarkably, a similar definition works with NCT, allowing us to define a ``Semi-Standard'' version of it. 

In particular, given a filling $\mu = (\mu_1, \mu_2, \ldots, \mu_k)$, consider segments $[1, \mu_1], [\mu_1+1, \mu_2], \ldots, [\mu_{k-1}+1, \mu_k]$ as above. Consider all NCT such that each of its parts contains at most one element from each segment. Define an equivalence relation on this set where two NCT are equivalent if one is obtained from the other by the action of the Young subgroup $S_{|\mu_1|} \times S_{|\mu_2|} \times \dotsc \times S_{|\mu_k|}$ which permutes points of the above segments. 

\begin{definition}
The {\it {Semi-Noncrossing Tableaux}}, or SNCT, is the set of equivalence classes we obtain in this way.
\end{definition}

\begin{theorem}
\label{semi}
The number of SSYT of shape $\lambda$ and weight $\mu$ is equal to the number of SNCT of shape $\lambda$ and weight $\mu$.
\end{theorem}

\begin{proof}
Given a Yamanuchi word $w$ and composition $\mu$, we say that $w$ SYT-agrees with $\mu$ if the SYT corresponding to $w$ agrees with $\mu$ as described above. Similarly define NCT-agreeing. Let $w_{\mu}$ be the largest word in lexicographic order obtained from $w$ by action of $S_{|\mu_1|} \times S_{|\mu_2|} \times \dotsc \times S_{|\mu_k|}$. Then the theorem follows from the following two statements, which are easy to verify: if $w$ is Yamanuchi and SYT-agrees with $\mu$, than $w_{\mu}$ is Yamanuchi; if $w$ is Yamanuchi and NCT-agrees with $\mu$, than $w_{\mu}$ is Yamanuchi. Now one can see that if $w_{\mu}$ is Yamanuchi then it SYT-agrees and NCT-agrees with $\mu$. Therefore both numbers in the statement of the Theorem are equal to the same number: the number of $w_{\mu}$-s which are Yamanuchi.
\end{proof}

We sometimes call SNCT just NCT if it is clear from context what is meant. It is possible to view SNCT as a filling of a Young Diagram. For example, here is a Figure showing a SNCT and the corresponding diagram, on which a potential $F$ is also marked: 

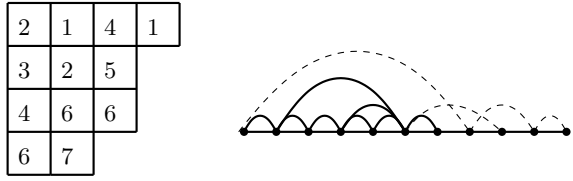
\begin{center}
\input{nct17.pstex_t}
\end{center}

\subsection{Relation to Temperley-Lieb Algebra} \label{TL}

Restrict our attention for a moment to partitions with two parts, that is $p=2$. Then our NCT are non-crossing matchings. The Garnir relations which generate all relations among $P_T$-s can then be described as follows. Given two crossing parts $a_1 < b_1 < a_2 < b_2$, we form two new Tableaux with parts $(a_1, b_1)$, $(a_2,b_2)$ and $(a_1, b_2)$, $(b_1, a_2)$. This corresponds to the equality: $$(a_1-a_2)(b_1-b_2) = (a_1-b_1)(a_2-b_2)+(a_1-b_2)(b_1-a_2).$$ We say that in this way we {\it {resolved}} the crossing $(a_1, a_2)$ and $(b_1, b_2)$. Given a matching, we can arbitrarily pick pair of crossing edges and resolve it. Note that the total number of crossings in each of the two resulting matchings is strictly smaller than in the original one. Thus after repeated application the process of resolving will stop. Observe that the following is a consequence of Theorem~\ref{thm:basis}.

\begin{corollary}
The result of the resolving of crossings via operation described above does not depend on the particular way it is done.
\end{corollary}

\begin{proof}
If $T_0$ is the original Tableaux there is a unique way to represent $P_{T_0}$ as a sum of $P_T$-s for NCT $T$. 
\end{proof}

There is a way to see directly that the resulting multiset of NCT is going to be the same no matter in which order we do the resolving. Given a Tableau $T$, take a crossing of two segments, ``\X'', and replace it with {\it vertical uncrossing\/} ``\noXv''and {\it horizontal uncrossing\/} ``\noXh''. When this is done to each crossing, we obtain the set $U_T$ of non-crossing matchings $u_i$, each possibly containing some cycles. For each of these matchings, let $T_i$ be the underlying NCT, and let $c_i = (-2)^{c(u_i)}$, where $c(u_i)$ is the number of cycles in $u_i$.

\begin{theorem}
$P_T = \sum_{i} c_i P_{T_i}$. 
\end{theorem}
 
\begin{proof}
 
It follows from the following ``Reidemeister moves'':

\begin{center}
\input{nct1.pstex_t}
\end{center}

\

\begin{center}
\input{nct2.pstex_t}
\end{center}

\end{proof}

Note that using this theorem we can compute the explicit entries of the images of elements of $S_n$ in an irreducible representation labeled by a two-part partition.

Recall that the {\it {Temperley-Lieb algebra}} $\TL_n(\xi)$ is the $\mathbb
C[\xi]$-algebra generated by $t_1, \ldots, t_{n-1}$
subject to the relations $t_i^2=\xi\, t_i$, and $t_i t_j t_i=t_i$ if
$|i-j|=1$, $t_i t_j = t_j t_i$ if $|i-j| \geq 2$. The dimension of
$\TL_n(\xi)$ equals the $n$-th Catalan number
$C_n=\frac{1}{n+1}\binom{2n}{n}$. 
A natural basis of the Temperley-Lieb algebra is 
$$\{t_w\mid w \textrm{ is a 321-avoiding permutation in } S_n\},$$
where $t_w := t_{i_1} \cdots t_{i_l}$, for a reduced decomposition $w = s_{i_1}
\cdots s_{i_l}$. In the case $q=1$ the map $\theta: s_i \mapsto t_i+1$ determines
a homomorphism $\theta:H_n(1) \to \TL_n(-2)$ from the Hecke algebra to the Temperley-Lieb algebra. The generators $t_i$ can be conveniently represented as {\it {Kauffman diagrams}}, shown in the Figure. 

\

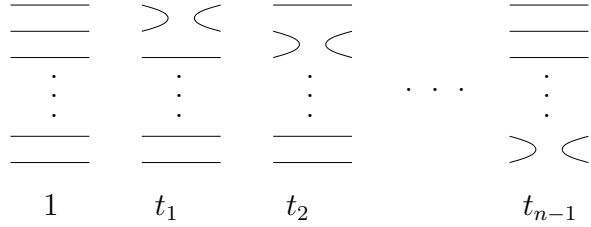
\begin{center}
\input{nct3.pstex_t}
\end{center}

The construction above is related to $\TL_n(-2)$. Indeed, assume we have a particular situation when edges of the Tableau $T$ form a matching between sets $1, \ldots, l$ and $l+1, \ldots, 2l$. Then $T$ can be viewed as a wiring diagram of some permutation $\omega(T) \in S_l$. There exists an obvious correspondence between NCT $T$ and generators $t_{u(T)}$ of $\TL_n(-2)$, since both correspond naturally to  non-crossing matchings on the set of $2l$ vertexes. Let $T_0$ be a NCT. 

\begin{theorem}
The coefficient of $P_{T_0}$ in the decomposition of $P_T$ is equal to the coefficient of $t_{u(T_0)}$ in $\theta(\omega(T))$.
\end{theorem}

\begin{proof}
One can see that the procedure of resolving crossings into two possible uncrossings corresponds exactly to choosing term $1$ or $t_i$ in the decomposition $\prod (t_i+1)$ of $\theta(\omega(T))$.
\end{proof}

Note that we thus obtain an alternative to Garnir relations, which are the tool for converting into SYT-basis.

\section{Standard and Non-Crossing Bitableaux} \label{sec:bi}
\subsection{Background on Standard Bitableaux}


Let $X = (x_{ij})$ be an $n \times n$ matrix. We consider minors $P_{I,J}$ of $X$ indexed by sets $I$ of rows and $J$ of columns. We consider monomials in those minors, $M = P_{I^1, J^1}, \cdots P_{I^h, J^h}$, called {\it {bideterminants}}. We say that the {\it content} of a monomial $M$ is $(\alpha, \beta)$ if the multiset of indices of rows in $M$ is $\alpha$, and the multiset of columns in $M$ is $\beta$. For example, if $M = \det (X)$ then $\alpha = \beta = \{1,2,\ldots,n\}$. While if $M = x_{12}^3$ then $\alpha = \{2,2,2\}$, $\beta = \{1,1,1\}$. Let $V(\alpha, \beta)$ be the vector space generated by monomials with content $(\alpha, \beta)$.

It is convenient to label monomials by {\it {bitableaux}}; see [Rota, et.al.]. Namely, a bitableau $(T, T')$ is a pair of Young Tableaux of the same shape $\lambda$ filled with positive integers (not larger than $n$). To obtain a monomial from a bitableau we take a column of $T$ and let its filling be the set $I$, and the filling of the corresponding column of $T'$ be $J$. Then $\prod P_{I,J}$ is the monomial corresponding to $(T,T')$. We can assume that entries in columns of $T$, $T'$ are distinct (otherwise the minor evaluates to $0$) and arranged increasingly.

A bitableau $(T, T')$ is {\it {standard}} if both $T$ and $T'$ are Semistandard Young  Tableaux, i.e. the numbers strictly increase in columns and weakly increase in rows. The following theorem is stated and proved in \cite{DKR}.

\begin{theorem}
The space $V(\alpha, \beta)$ is generated by monomials corresponding to standard bitableaux of content $(\alpha, \beta)$
\end{theorem}

The following is the simplest example of this fact: $$x_{12}x_{21} = x_{11}x_{22} - (x_{11}x_{22} - x_{12}x_{21})$$ which can be schematically written as shown in the following Figure:

\begin{center}
\input{nct14.pstex_t}
\end{center}

\subsection{Non-Crossing Bitableaux}

In this Section it is more convenient to use the modified definition of NCT with pre-fixed $F$ which was given in previouse Section. 

Call a bitableau $(T, T')$ non-crossing if both $T$ and $T'$ are non-crossing.

The following Figure is an example of non-crossing bitableau of content $$(\{1,1,2,2,3,4,4,5,6,6,6,7\},\{1,1,2,2,2,3,4,5,5,5,6,7\})$$:

\begin{center}
\input{nct15.pstex_t}
\end{center}

We are ready to state the main theorem of this section.

\begin{theorem}
\label{bi}
The space $V(\alpha, \beta)$ is generated by monomials corresponding to non-crossing bitableaux of content $(\alpha, \beta)$.
\end{theorem}

In order to prove it we first prove two lemmas.

\begin{lemma} \label{le}
The number of non-crossing bitableaux of given content $(\alpha, \beta)$ and shape $\lambda$ is equal to the number of standard bitableaux of the same content and shape.
\end{lemma}

\begin{proof}
Follows from the Theorem \ref{semi} and the fact that it is possible to find appropriate $F$. 
\end{proof}

In fact, the objects are in bijection. For example, here is the (semi) standard bitableau corresponding to the non-crossing tableau  above.

\begin{center}
\input{nct16.pstex_t}
\end{center}

\begin{lemma} \label{lu}
It is enough to prove Theorem \ref{bi} for the case when $\alpha$ and $\beta$ do not contain repeated elements.
\end{lemma}

\begin{proof}
Assume $\alpha$ and $\beta$ are not multiplicity free. Add some rows and columns to  matrix $X$ so that now they are. Call such an operation {\it {cloning}}. For example, the bitableau in the Figure above is defined on a matrix with (at least) $7$ rows and at least $7$ columns. The cloning procedure would consist of constructing a new matrix by taking column $1$ twice, column $2$ twice, column $4$ twice, and column $6$ - thrice. Similarly row $1$ twice, row $2$ and row $5$ thrice (compare with the content of this bitableau mentioned above). The new matrix has the property that the same bideterminant can be now written without using the same column or the same row of the matrix twice. 

Then, assuming Theorem \ref{bi} holds for $\alpha$ and $\beta$ multiplicity free, we can decompose any monomial $M$ into non-crossing monomials. The monomials which contain a pair of cloned rows or columns are equal to $0$. The rest is going to remain non-crossing when we go back to original pre-cloning matrix. Thus non-crossing bitableaux span $V(\alpha, \beta)$. However, the number of them is equal to the number of standard bitableaux as asserted by Lemma \ref{le}. Since this number is the dimension of $V(\alpha, \beta)$ as a vector space, we conclude the Theorem \ref{bi} is true for $(\alpha, \beta)$ which are not multiplicity-free.
\end{proof}

Now we are ready to prove Theorem \ref{bi}.

\begin{proof}
According to Lemma \ref{lu} it is enough to prove the statement for the case of multiplicity-free $(\alpha, \beta)$. The number of non-crossing bitableaux equals the dimension of $V(\alpha, \beta)$, as we know from the Lemma \ref{le}. Thus it is sufficient to show that monomials corresponding to non-crossing bitableaux are linearly independent. 

Let $a_{i1} < \cdots < a_{ik_i}$ be the $i$-th column of a non-crossing tableau $T$, $b_{i1} < \cdots < b_{ik_i}$ be the $i$-th column of a non-crossing tableau $T'$. Restrict our attention to the matrices $X$ such that for each $j$ the columns of $X$ indexed by $a_{ij}$ are equal for all $i$ and the rows of $X$ indexed by $b_{ij}$ are equal for all $i$. Obviously, it is sufficient to show linear independence for such matrices. 

Introduce the following {\it {lexicographic}} order on non-crossing Tableaux: take the corresponding readings as in Section \ref{NCT} and order them lexicographically. The key observation is that if we restrict our attention to matrices corresponding to $(T,T')$ as described above, then any non-crossing bitableau $(R, R')$ takes non-zero value on those matrices only if $R \leq T$, $R' \leq T'$ in that order. This is exactly the same argument as in the proof of Theorem \ref{thm:basis}. 

Pick {\it {the smallest}} in the above lexicographic order non-crossing bitableau $(T,T')$ (we can assume that we order first by $T$, then by $T'$) such that there is a linear relation among non-crossing bitableaux involving this particular monomial with a non-zero coefficient. Restrict attention to matrices corresponding to $(T,T')$ as above. Since the monomial corresponding to $(T,T')$ does not evaluate to $0$ identically on these matrices, there should be another monomial $(R, R')$ in this linear relation with the same property. The only possibility for $(T,T')$ to be minimal while satisfying properties $R \leq T$, $R' \leq T'$ is $(T,T')=(R, R')$, which is a contradiction. Thus linear dependences are not possible and the proof is complete.
\end{proof}

\section{Standard and Non-Crossing Monomials} \label{sec:mon}

Let us recall the Standard Monomial Theory for the coordinate ring of the Grassmanian, closely following the way it is presented in \cite{M}.

Denote by $\mathbb G_{m,n}$ the set of linear subspaces of dimension $m$ of $\mathbb C^{m+n}$. The set $\mathbb G_{m,n}$ is called (complex) {\it {Grassmannian}}. Once a basis for $\mathbb C^{m+n}$ is chosen, $\mathbb G_{m,n}$ can be identified with the set of $m \times (m+n)$ matrices of rank $m$. For a set $I = \{i_1, \ldots, i_m\}$, let $P_I$ be the determinant of the minor of an $m \times (m+n)$ matrix which corresponds to columns $i_1, \ldots, i_m$. The rank condition on matrices in $\mathbb G_{m,n}$ is equivalent to the polynomial system $\forall I \; P_I \not = 0$. Thus $\mathbb G_{m,n}$ is a complex algebraic variety. The $P_I$-s are called {\it {Plucker coordinates}}. 

Fix a complete flag $0 = V_0 \subset V_1 \subset \cdots \subset V_{m+n} = \mathbb C^{m+n}$. Let $\lambda$ be a partition contained in an $m \times n$ rectangle. The {\it {Schubert variety}} $X_{\lambda}$ is defined as follows: $$X_{\lambda} = \{W \in \mathbb G_{m,n} \mid \dim(W \cap V_{n+i-\lambda_i}) \geq i, \; 1 \leq i \leq m\}.$$

Let $\{i_1, \ldots, i_m\}$ and $\{j_1, \ldots, j_m\}$ be two sets of integers between $1$ and $m+n$; let $l$ be an integer between $1$ and $m$. Let $S$. $S'$ and $S''$ be the groups of permutations of $i_l, \ldots, i_m, j_1, \ldots, j_l$, $i_l, \ldots, i_m$ and $j_1, \ldots, j_l$ respectively. Let $S/(S' \times S'')$ be the set of coset representatives of the Young subgroup $S' \times S''$. Finally, let $\epsilon(\omega)$ be the sign of permutation $\omega$. Then on $\mathbb G_{m,n}$ we have $$\sum_{\omega \in S/(S' \times S'')} \epsilon(\omega) P_{i_1, \ldots, i_{l-1}, \omega(i_l), \ldots, \omega(i_m)} P_{\omega(j_1), \ldots, \omega(j_m), j_{l+1}, \ldots, j_m} = 0.$$ These are the so called {\it {Plucker relations}}. 

\begin{theorem}
{\rm [Theorem~3.1.6, \cite{M}]{}} \
The Plucker relations completely determine the Grassmanian and generate the ideal $I(\mathbb G_{m,n})$ of the variety.
\end{theorem}

Define a lexicographic order on $m$-tuples $J = j_1 < \cdots < j_m$ of integers: $J < J'$ iff $\exists \; k$ s.t. $j_k<j'_k$, $j_l = j'_l$ for $l<k$. Let $M = P_{J^1}, \cdots P_{J_h}$ be a monomial in Plucker coordinates. We represent $M$ by the Tableaux 

$$
\begin{array}{ccc}
j_1^1  & \ldots & j_m^1 \\
\ldots & \ldots & \ldots \\ 
j_1^h & \ldots & j_m^h 
\end{array}
$$

We call $M$ a {\it {standard monomial}} if $J^1 \leq \cdots \leq J^h$. 

The following theorem reveals the reason for defining standard monomials. Let $R_{m,n} = \mathbb C[P_J]/I_{m,n}$, where $I_{m,n}$ is the ideal generated by Plucker relations.

\begin{theorem} \label{smcr}
{\rm [Theorem~3.3.4, \cite{M}]{}} \
The ideal $I_{m,n}$ is radical, thus $R_{m,n}$ is the coordinate ring of the Grassmannian. Standard monomials form a basis of $R_{m,n}$.
\end{theorem}
 
In particular, any monomial in Plucker coordinates can be (uniquely) expressed as a linear combination of standard monomials. 

The following theorem is a stronger statement. Let $\lambda \subset n \times n$ be a partition. Define the $m$-tuple $I$ by $i_k=n+k-\lambda_k$. 

\begin{theorem} \label{SV}
{\rm [Theorem~3.3.4, \cite{M}]{}} \
Standard monomials $M$ s.t. $M P_I$ is also standard form a basis for the coordinate ring $R_{\lambda}$ of the Schubert variety $X_{\lambda}$.
\end{theorem}

Call a monomial $M$ {\it {non-crossing}} if the following condition fails to hold for any $p,q,k$: $j_k^p<j_k^q<j_{k+1}^p<j_{k+1}^q$. In other words, $M = P_{J^1}, \cdots P_{J^h}$ is non-crossing if sets $J^k$ are non-crossing, in the terminology above. 

\begin{theorem}
\label{plu}
Every monomial in Plucker coordinates is a linear combination of NCM, and in a unique way. NCM form a basis for the coordinate ring of the Grassmannian.
\end{theorem}

\begin{proof}
Non-Crossing Monomials are a particular case of bideterminants of Section \ref{sec:bi}. Indeed, they are the bideterminants corresponding to pairs $(T,T')$ where $T$ is a square Tableaux with row $i$ filled with $i$-s. Then Lemma \ref{le} implies that the number of NCM is equal to the dimension of $R_{m,n}$. Theorem \ref{bi} implies that NCM-s are linearly independent. We conclude that they form a basis for $R_{m,n}$.
\end{proof}

An analog of Theorem \ref{SV} can also be stated. Call $J = (j_1< \cdots < j_n)$ {\it {relevant to $\lambda$}} if for all $k$ $j_k \leq n+k-\lambda_k$. Call a monomial $M = \prod P_{J^i}$ {\it {relevant to $\lambda$}} if each $P_{J^i}$ is relevant to $\lambda$. Note that a Standard Monomial $M$ is relevant to $\lambda$ iff $M P_I$ is also a Standard Monomial, where $I$ is constructed from $\lambda$ as described above.

\begin{lemma} \label{rel}
The number of NCM relevant to $\lambda$ is equal to the number of Standard Monomials relevant to $\lambda$.
\end{lemma}

\begin{proof}
It is easy to see that a bijection between SYT and NCT described after Corollary \ref{eq} does not change the relevance to $\lambda$.
\end{proof}

\begin{theorem}
NCM relevant to $\lambda$ form a basis for the coordinate ring $R_{\lambda}$ of the Schubert variety $X_{\lambda}$.
\end{theorem}

\begin{proof}
The dependence relations that hold in $R_{m,n}$ must also hold in $R_{\lambda}$. Thus Theorem \ref{plu} implies that every element in $R_{\lambda}$ is a linear combination of NCM. Note that all monomials not relevant to $\lambda$ vanish in $R_{\lambda}$. Therefore every element in $R_{\lambda}$ is a linear combination of NCM relevant to $\lambda$.

Theorem \ref{SV} and Lemma \ref{rel} imply that the number of NCM relevant to $\lambda$ is equal to the dimension of $R_{\lambda}$ as a vector space. Thus, we conclude that NCM relevant to $\lambda$ form a basis for $R_{\lambda}$.
\end{proof}

\subsection{Gr\"{o}bner basis}

The idea of this subsection was suggested by David Speyer. We follow \cite{St} and Chapter 14 of \cite{MS}. 

A monomial $\bf x^{\bf a} = x_1^{a_1}\cdots x_n^{a_n}$  in $k[\bf x]$ corresponds to the point $\bf a \in \mathbb N^n$. A monomial order is a total order on $\mathbb N^n$ such that $(0, \ldots, 0)$ is a minimal element and $\bf a < \bf b$ implies $\bf {a+c} < \bf {b+c}$ for any $\bf c$. An {\it {initial}} monomial, or term $in_<(f)$ of a polynomial $f$ is a minimal monomial with non-zero coefficient in $f$ with respect to given monomial order $<$. If $I$ is an ideal, the {\it {initial ideal}} $in_<(I)$ is the ideal generated by initial monomials of all elements of $I$. The {\it {Gr\"obner basis}} is a subset $G \subset I$ such that $in_<(g), g\in G$ form a basis for $in_<(I)$. If no monomial in $in_<(g), g\in G$ is redundant, $G$ is said to be minimal. Non-initial monomials in $I$ are called {\it {standard}}, we prefer to call them just non-initial though to avoid ambiguity in notation.

\begin{proposition}\label{Gr}
{\rm [Proposition~1.1, \cite{St}]{}} \
The (images of) the non-initial monomials form a basis for $k[\bf x]/I$
\end{proposition}

This can be used to show that Standard monomials form a basis for $R_{m,n}$. Indeed, introduce the lexicographic order $<$ on tuples $J$ indexing Plucker coordinates $P_J$. Let $<$ also denote reverse lexicographic term order induced by $<$. 

\begin{theorem} \label{sg}
{\rm [cf. Theorem~14.6, \cite{MS}]{}} \
The products $P_J P_{J'}$ for nesting $J$ and $J'$ generate the initial ideal $in_<(I_{m,n})$ of the ideal of Plucker relations with respect to $<$.
\end{theorem}

\begin{remark}
Note that Theorem~14.6, \cite{MS} is stated and proved in the more general setting of {\it {Plucker algebra}}.
\end{remark}

Theorem \ref{sg} together with Proposition \ref{Gr} imply Theorem \ref{smcr}. It is natural to ask whether Theorem \ref{plu} can be derived in this way. This poses the following 

\begin{question} \label{qg}
Is there monomial order such that NCM are exactly the non-initial terms of the ideal of Plucker relations?
\end{question}

We answer this affirmatively for the case of $\mathbb G_{2,n}$.

Let $J=(j_1, j_2)$, and denote weight $w(J) = j_2-j_1$. Let $w(P_{J^1} \cdots P_{J^h}) = w(J^1) \cdots w(J^h)$. Order monomials reversely with respect to their weight, that is the larger weight the smaller is the monomial. Let order be arbitrary in case of equal weight. Denote this order by $\prec$. 

\begin{theorem}
NCM are exactly the non-initial monomials with respect to $\prec$.
\end{theorem}

\begin{proof}
First we show that crossing (e.i. not non-crossing) monomials lie in $in_{\prec}(I_{2,n})$. Indeed, assume $J$ and $J'$ in decomposition $M = \prod P_{J^k}$ do cross. Without loss of generality we can interpret that as $j_1 < j'_1 < j_2 < j'_2$. Then in relation $P_{(j_1, j'_1)} P_{(j_2, j'_2)} - P_{(j_1, j_2)} P_{(j'_1, j'_2)} + P_{(j_1, j'_2)} P_{(j'_1, j_2)}$ the term $P_{(j_1, j_2)} P_{(j'_1, j'_2)}$ has the largest weight $(j_2-j_1)(j'_2-j'_1)$ and this is the smallest with respect to $\prec$. Then $P_J P_{J'}$ lies in $in_{\prec}(I_{2,n})$, and thus $M$ also lies in $in_{\prec}(I_{2,n})$.

It remains to show that $in_{\prec}(I_{2,n})$ has no other monomial generators. However we already know that, since we know that NCM are linearly independent. This finishes the proof.
\end{proof}

\section{Non-Crossing Tableaux and Irreducible Representations of $GL_n$} \label{GL}

Let $X = (x_{ij})$ be an $n \times r$ matrix. Given $\lambda$ with at most $n$ rows, let $T_{\lambda}$ be a Tableau of shape $\lambda$ with row $i$ filled with $i$-s. Let $T$ be any Tableau of shape $\lambda$ and let $M_{T}$ be the bideterminant corresponding to $(T_{\lambda}, T)$ as in Section \ref{sec:bi}. Let $D_{\lambda}$ be the space spanned by all such $M_T$. Note that $GL_n(\mathbb C)$ acts on $X$ by left multiplication, which makes $D_{\lambda}$ into a $GL_n$-module. 

The following is a reformulation of Theorems 4.5 and 4.7b in \cite{Gr}. See also Theorem 2 and Corollary after Theorem 1 in Chapter 8, \cite{F}.

\begin{theorem}
$D_{\lambda}$ is the irreducible $GL_n$-module corresponding to highest weight $\lambda$. The $M_T$-s for Semi-Standard Young Tableaux $T$ form a basis for $D_{\lambda}$.
\end{theorem}

Now, define Non-Crossing Tableaux as in Section \ref{sec:bi} (these are Semi-Non-Crossing Tableaux in the terminology of Section \ref{NCT}).

\begin{theorem}
The $M_T$-s for Non-Crossing Tableaux $T$ form a basis for $D_{\lambda}$.
\end{theorem}

\begin{proof}
Lemma \ref{le} implies that the number of NCT is equal to the dimension of $D_{\lambda}$. Theorem \ref{bi} implies that $M_T$-s are linearly independent. Combining we conclude that they form a basis.
\end{proof}

\section{Concluding Remarks} \label{goodbye}

Observe that all the proofs in the paper were very much non-constructive. In particular, unlike the case of standard objects, we do not have a constructive way to write any given element of the appropriate space as a linear combination of non-crossing basis elements. The only exception comes from the Temperley-Lieb algebra described in Subsection \ref{TL}. 

\begin{question}
Can the results of Subection \ref{TL} be generalized to the case of Tableaux with more than $2$ rows? 
\end{question}

\begin{question}
Is there an uncrossing algorithm in the spirit of the straightening algorithm of Desarmenien \cite{D} which decomposes a given Bitableau into Non-Crossing Bitableaux?
\end{question}

An affirmative answer to Question \ref{qg} would implicitely provide such an algorithm.

The simplest instance of this question be stated as follows. Consider a $3 \times 3n$ matrix. Split the columns into $n$ parts with $3$ columns in each. Pictorially we represent it by a diagram consisting of $3n$ dots on a line, with each $3$-set represented by $2$ arcs connecting the leftmost and the rightmost dots with the middle one. We call such a triple a {\it {seagull}}, and we call the arcs {\it {left wing}} and {\it {right wing}}. To each seagull corresponds a $3 \times 3$ minor of the matrix, and to a particular partition into seagulls corresponds the product of corresponding minors.

We are allowed to do the following procedure: pick two seagulls among the whole set, and using some Plucker relation substitute them for a linear combination of (three) other pairs of seagulls. One example can be found in Figure \ref{nct7}. In that figure one can see that the resulting pairs of seagulls are non-crossing. It follows for example from Theorem \ref{plu} that it is always possible to find an appropriate series of such moves so that all partitions in the resulting linear combination are non-crossing. 

\begin{question}
What is the algorithm for finding such a decomposition, and why does it terminate?
\end{question}

Note that unlike the case of similar question for $2 \times 2$ minors of a $2 \times 2n$ matrix, it is not necessarily true that the number of crossings can be reduced at each step! 

The author contacted a number of people with questions regarding this work. He is grateful to the following people for making comments that influenced the development of the paper: Alexander Postnikov, Arun Ram, David Speyer, Thomas Lam, Richard Stanley, Andrew Mathas. The author is also very grateful to Denis Chebikin who helped with editing of the paper.

\end{document}